\documentclass{llncs}
\usepackage{subcaption}
\usepackage[pdftex]{graphicx}
\usepackage[usenames,dvipsnames]{color}
\usepackage{overpic}
\usepackage{wrapfig}
\usepackage{svg}
\usepackage{hyperref}
\usepackage[normalem]{ulem}
\usepackage[superscript,biblabel]{cite}
\newcommand{\figref}[1]{Fig.~\ref{#1}}
\newcommand{\secref}[1]{Section~\ref{#1}}
\newcommand{\blender}{\emph{Blender}}
\newcommand{\sofa}{\emph{SOFA}}
\newcommand{\vid}{(\jorg{see Video, Supplemental Digital Content 1})}
\newenvironment{packed_item}{
\begin{itemize}
\setlength{\itemsep}{0pt}
\setlength{\parskip}{0pt} 
\setlength{\parsep}{0pt}
}{\end{itemize}}
\newcommand{\jorg}[1]{\textcolor{blue}{#1}}

\begin{document}

\title{Adding Safety Rules to Surgeon-Authored VR Training}

\titlerunning{TIPS safety assessment}  

\author{
Ruiliang Gao\inst{1}, Sergei Kurenov\inst{2},  Erik W. Black\inst{1,3} \and J\"org Peters\inst{1,3,4}}
%
%
\institute{$^1$ University of Florida, Gainesville, FL32611, USA,\\
$^2$ Roswell Park Comprehensive Cancer Center,
Buffalo, NY14263, USA\\
$^3$ Ph.D., Professor\\
$^4$ corresponding author: jorg.peters@gmail.com\\
432 Newell Dr, U of Florida,
Gainesville, FL, 32611-6120,USA,\\ (352)382 1200 (fax 1220) }


\maketitle              


\begin{abstract}
\textcolor{white}{.}
\\
\textbf{Introduction}
Safety criteria in surgical VR training are typically hard-coded and informally summarized.
The Virtual Reality (VR) content creation interface, TIPS-author, for the 
Toolkit for Illustration of Procedures in Surgery (TIPS) allows surgeon-educators
(SE s) to create laparoscopic VR-training modules with force feedback.
TIPS-author initializes anatomy shape and physical properties selected by the SE accessing a cloud data base of physics-enabled pieces of anatomy. 
\\
\textbf{Methods}
A new addition to TIPS-author are safety rules
that are set by the SE and are
automatically monitored during simulation. Errors are recorded as visual snapshots for feedback to the trainee.
This paper reports on the implementation and opportunistic evaluation of the snapshot mechanism
as a trainee feedback mechanism.
%
TIPS was field tested at two surgical conferences, one before and one after adding the snapshot feature.
\\
\textbf{Results}
While other ratings of TIPS remained unchanged for an overall Likert scale score of 5.24 out of 7 (7 = very useful), the 
the rating of
`The TIPS interface helps learners understand the force necessary to explore the anatomy' improved from 5.04 to 5.35 out of 7 after the snapshot mechanism was added.
\\
\textbf{Conclusions}
The ratings indicate the viability of the TIPS open-source SE-authored surgical training units.
Presenting SE-determined 
procedural missteps via the snapshot mechanism at the end of the training increases acceptance.



\end{abstract}
\section{Introduction}
Teaching laparoscopic surgery under one-on-one supervision in the operating room (OR) is costly ranging, already a decade ago, from \$50--\$135 per minute \cite{Macario:2010}.
Less supervision is risky:
cauterizing too close to a sensitive organ or nicking a central vein are difficult to repair and may cause the patient unnecessary suffering.
Therefore alternative training methods are ethically and fiscally prudent.
Mentored self-study curricula, such as
Fundamentals of Laparoscopic Surgery (FLS), offer
dexterity training and certification on peg-board transfer, cutting and suturing  of physical props as a foundation before working on real patients \cite{Zendejas:2013}. 
However, FLS box training can not prepare for the high variability of anatomy and soft tissue response that actual cases present and
provides no automatic checking of safety criteria.

The additional technical challenge that this work addresses is that
entry must interpret and trigger deployment of monitors for a palette of surgical safety criteria set by surgeon educators.
To explain the challenge,
we first review  soft tissue simulation and VR trainers in general in \secref{sec:review}, then a customizable training framework (TIPS, \secref{sec:TIPS}) and then, in \secref{sec:missing}, formulate
the specific challenge and new contributions.

\subsection{Soft tissue simulation and VR trainers}
\label{sec:review}

The last decade has witnessed progress in soft tissue simulation for a range of surgical scenarios such as laparoscopic surgery, heart surgery, neurosurgery, orthopedic and arthroscopic
surgery.
Early multilayered tissue models
for orthopedic trauma surgery
were based on 3D mass-spring systems accelerated with graphics hardware
\cite{qin2010novel}.
More recently, simulation of cardiac electrophysiology simulation, pre-operative planning of cryosurgery and per-operative guidance for laparoscopy use finite elements in real time
in the open source SOFA soft tissue simulation
platform \cite{talbot2015surgery}.
A real-time neurosurgery simulator of skull drilling and surgical interaction with the brain was proposed in \cite{echegaray2014brain}
and
\cite{arikatla2018high} report on simulation of 
drilling and cutting of the bone using the burr and the motorized oscillating saw based on the open source iMSTK framework.
\cite{mitchell2015gridiron} presented a framework for interactive outlining of regions for 
simulation of reconstructive plastic surgery and
\cite{cecil2018advanced} describes a virtual surgical environment for training residents in less invasive stabilization system surgery used to address fractures of the femur.

Several commercial VR training environments aim to reduce time spent teaching in the OR by offering
training modules with virtual anatomy that can be probed using force feedback devices. 
Manual laparoscopic techniques lend themselves particularly well to simulation that leverage force-feedback devices.
Virtual reality simulators allow trainees to practice decision-making
and execution prior to entering the OR \cite{Satava2005,Gurusamy2009}.
A number of commercial solutions have sunset during the past 20 years (e.g.\ Simsurgery), or were merged or bought up by larger companies (see e.g.\ SurgicalScience, Simbionics, Mimic). However, 
commercial training environments neither capture the broad spectrum of physical variations encountered in laparoscopic practice, nor prepare learners for less common  
interventions.

\def\wid{1.0\linewidth}
\begin{figure}[ht]
\begin{subfigure}{0.55\textwidth}
\centering        
\includegraphics[width=\wid]{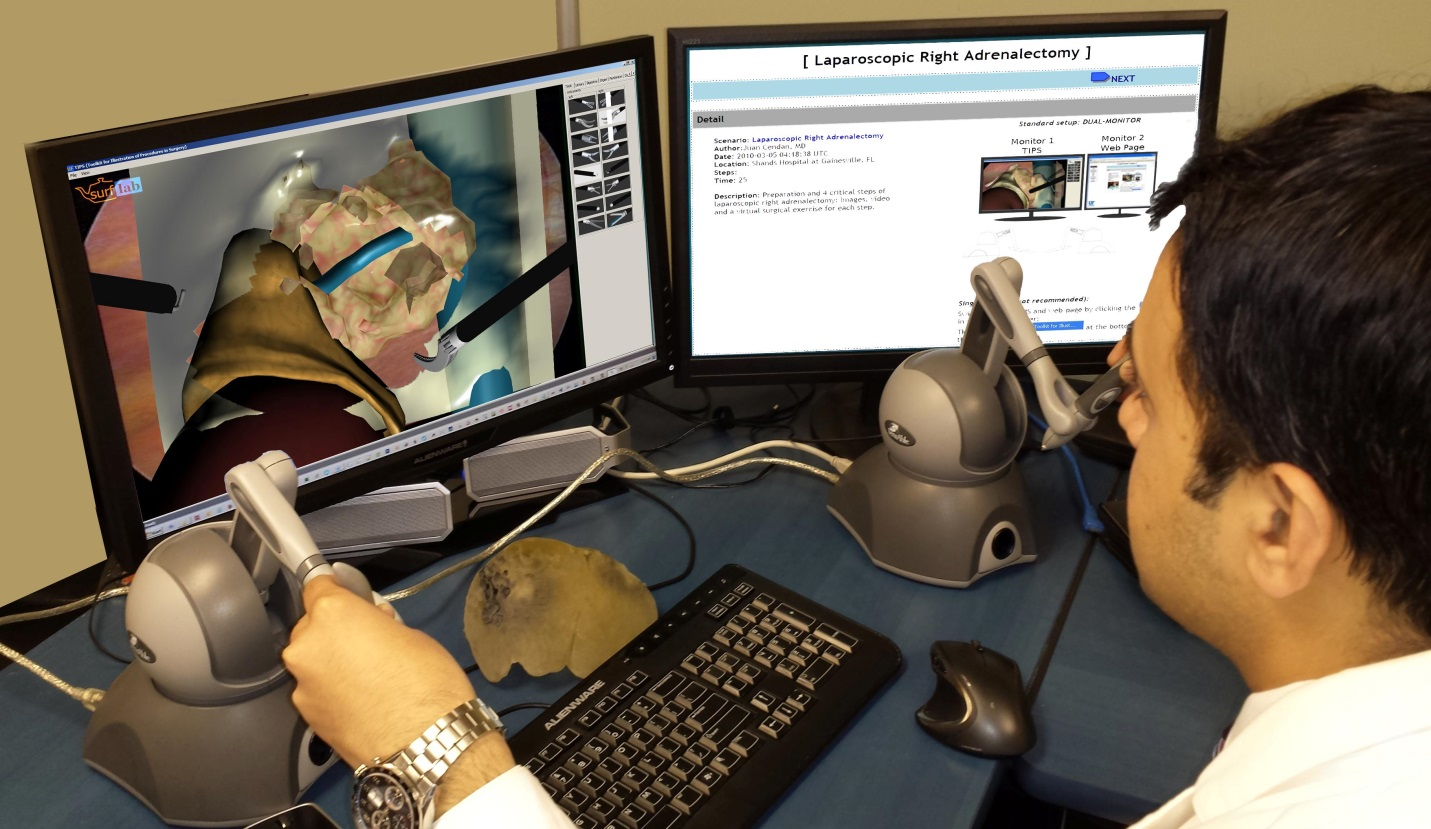}
\caption{
TIPS environment, no hand piece
}
\end{subfigure}
\hskip 0.05\linewidth
\begin{subfigure}{0.42\textwidth}
\centering        
\includegraphics[width=\wid]{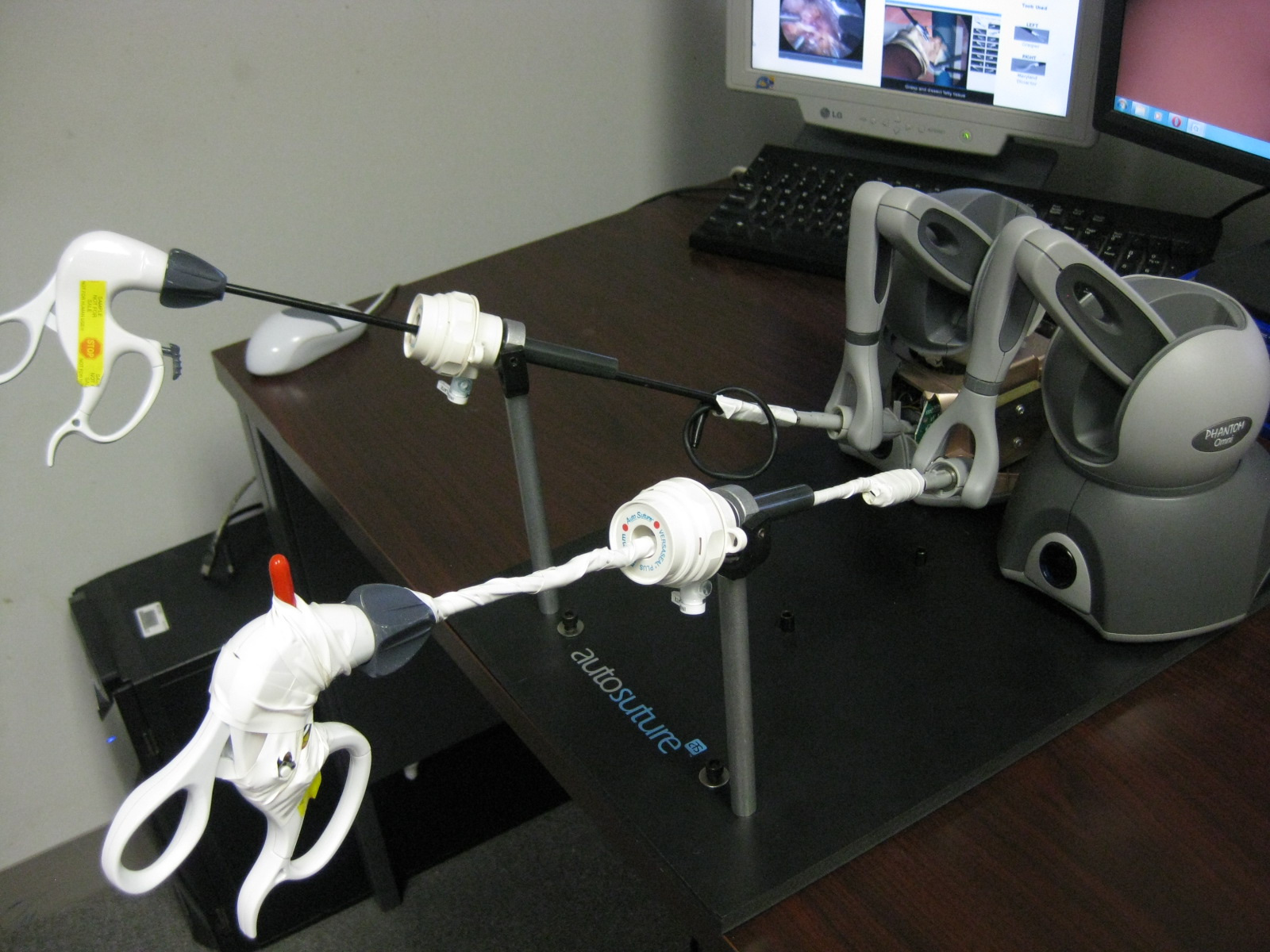}
\caption{setup with wired handles }
\end{subfigure}
\caption{Physical setup of the TIPS environment: (a) {\it left screen} TIPS-simulator; {\it right screen} TIPS-trainee
providing instructions and quizzes.
Two six degrees-of-freedom haptic devices provide physical feedback.
The haptic setup can  be augmented by  (b) laparoscopic
tool handles for a training setup with fulcrum, similar to an FLS box trainer, but wiring is cumbersome and brittle.
}
\label{fig:TIPSsetup}
\end{figure}

Building a robust virtual environment is  a formidable challenge,
leveraging scientific advances in collision detection, real-time differential equation solving, interactive visual and haptic feedback in a well-engineered interface. Creating or updating VR scenarios
is a second hurdle: due to the back-and-forth between engineers, computer scientists and medical experts, content creation is neither cheap nor fast and can take months or years; and the result is not easily adjustable to create variants, uncommon or specialized scenarios.

\def\wid{0.95\linewidth}
\def\widd{0.45\linewidth}
\subsection{The Toolkit for Illustration of Procedures in Surgery (TIPS)}
\label{sec:TIPS}
\begin{wrapfigure}{r}{0.65\linewidth}
\centering
\vskip-0.7cm
\includegraphics[width=\wid]{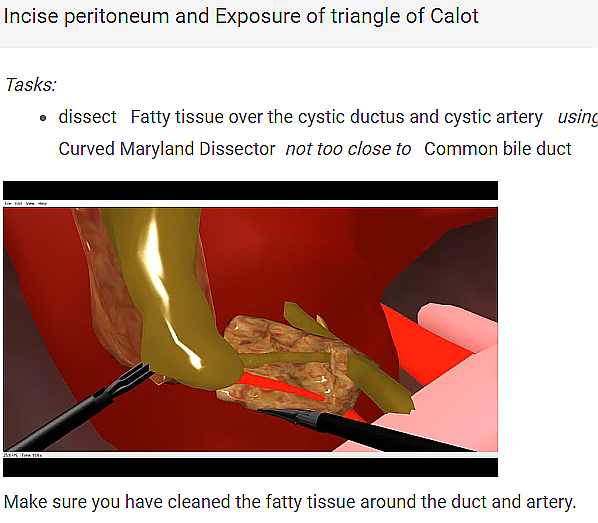}
\caption{Screenshot of TIPS-trainee instructions directly generated from the TIPS-author entries with embedded
(real and) simulation footage video clips.}
\label{fig:trainee}
\end{wrapfigure}
The open source Toolkit for Illustration of procedures in Surgery (TIPS) addresses the fast prototyping challenge of missing variants of anatomy and of less common
laparoscopic procedures.
\figref{fig:TIPSsetup}
shows two physical setup to transmit force via small robotic arms.
The TIPS open source  environment 
 consists of TIPS-simulator, an interactive soft-tissue laparoscopic simulation
with force feedback, (TIPS-trainee, a web-based component providing instruction and examples to a novice surgeon and
TIPS-author (see \secref{sec:author})  that allows the surgeon-educator to specify steps of a minimally invasive procedure in the fixed format:

The triple `action, anatomy, tool' is used to initialize the geometry and physical properties of the virtual  anatomy.
This high-level initialization is possible thanks to a rich database of simlets.
A simlet is a piece of anatomy with its physical properties, created by content artists. Simlets combine in a Lego-like fashion  \cite{Yeo:2010:MMVR,Dindar:2016:MMVR,Sarov:2018:AIS}. 
For example, the cystic duct, cystic artery and the fatty tissue covering them (each with their unique Young's modulus etc.), form an anatomy simlet.
In our implementation, the web-based TIPS-author  interface 
auto-completes typed items once they are recognized to be in the database and thereby steers the author towards existing simlets 
\vid.
Both the listing of steps and the resulting simulation is peer-reviewed for completeness and relevance before roll-out to trainees.


An initial study of 34 medical students \cite{Lesch:2019:VR:}
assessed whether the interactive learning within the prototype TIPS environment has advantages over passive learning from professional instructional videos \cite{wise}.
The study showed that inter-active TIPS platform
instilled greater confidence in the ability to reproduce the steps of the procedure  (p=0.001) and was preferred by the participants as a learning tool (p=0.011).
Of course confidence is not always positively correlated with
proficiency \cite{McGrane:2016}.

\subsection{The missing component and contribution: automatic initialization and monitoring of safety criteria}
\label{sec:missing}

The \emph{missing component in earlier work} is a lack of automatic interpretation of a `safety' entry in the 
specification of a surgical task. 
The technical challenge is that this entry must interpret and trigger deployment of monitors for a palette of surgical safety criteria set by surgeon educators. Examples are: do not cauterize near sensitive organs, 
limit the force when separating vessels from fatty tissue, etc.\ . The challenge addressed in this paper is to provide a \emph{generic
mechanism} for meaningful author-controlled yet automatic (unsupervised) surgical safety feedback to the trainee  -- and so accommodate current and future
not yet specified laparoscopic training scenarios.

While on one hand, haptic interactive simulation is no different from
(a specialized) computer game, on the other hand  the safety criteria can not be hard-coded a priori 
as in a game environment -- because the surgeon-author sets the criteria.
The challenge then is to insure that (a large class) of surgical safety criteria
can be 
(i) formulated, 
(ii) automatically translated into measurable events during VR 
simulation and 
(iii) generate both immediate feedback to the trainee and a meaningful final report to be shared with the instructor.
Section 2 introduces TIPS.
Section 3 reports on the \emph{new contributions} that
\begin{packed_item}
\item[1.]
allow the surgeon-educator to specify the surgical safety criteria;
\item[2.]
 are automatically monitored within TIPS-simulator;
\item[3.]
provide immediate feedback to the trainee; and
\item[4.]
return, in a secure fashion, visual feedback to both the trainee and
a summary message to the instructor
as a series of snapshots. 
\end{packed_item}

The list of errors and task-completion in the message (4.) record progress: a repeatedly empty error-report and completed task report indicate proficiency with respect to the training unit. This can be used to trigger the final assessment by the instructor and complete a feedback loop to improve the teaching unit by setting additional safety criteria or better specifying steps of the procedure.


\begin{figure}[ht]
\begin{center}
\includegraphics[width=1.0\textwidth]{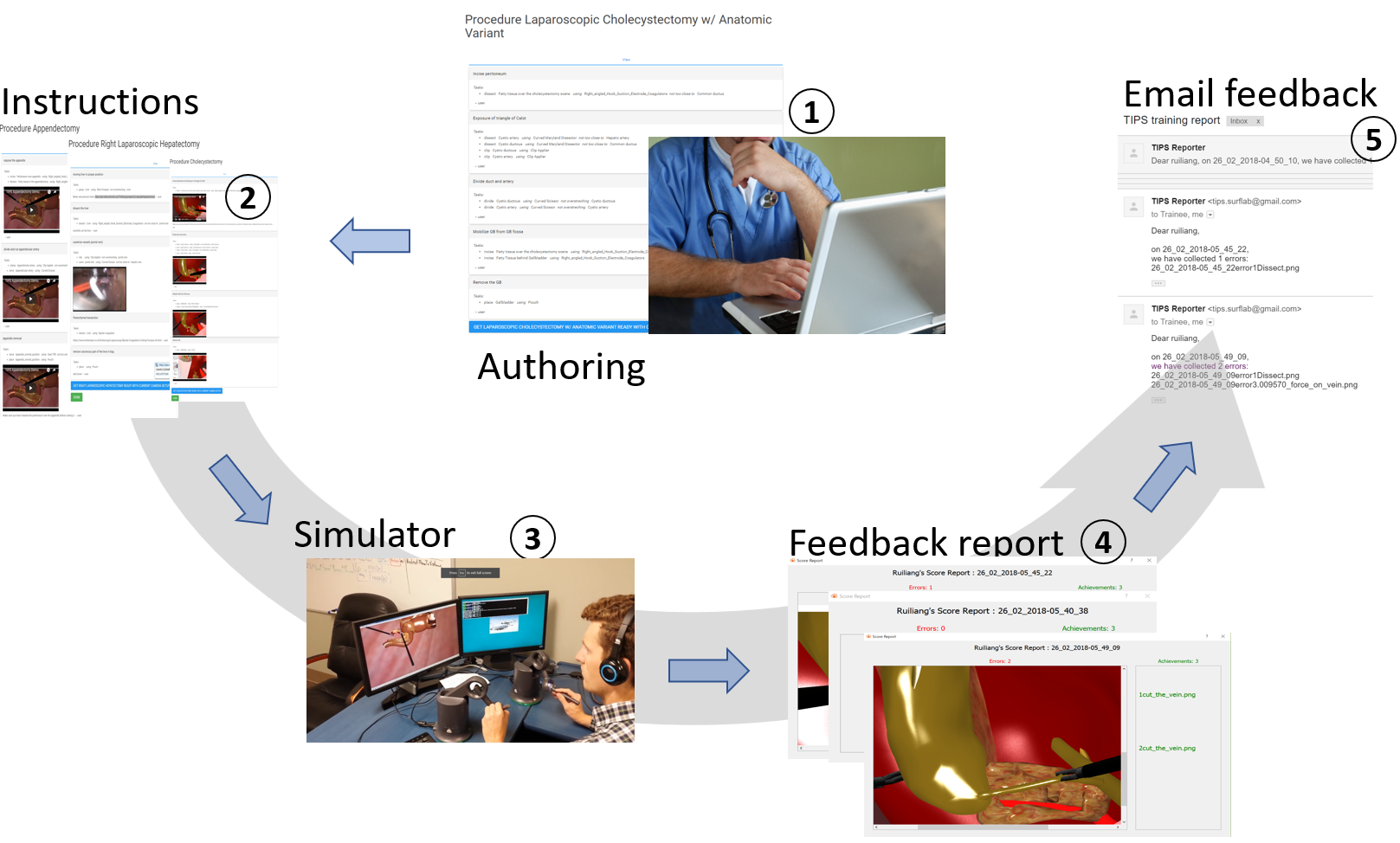}
\caption{TIPS-author defines the interactive VR  training simulation: (1) The author specifies procedural steps and safety concerns in a fixed format. (2) Instruction pages are generated from the author's description. (3) Simlets (pieces of anatomy and their physical properties) are combined to initialize the scenario in TIPS-simulator. (4) Trainee achievements and safety violations are screen-captured in TIPS-simulator for post-review.
This is the focus and contribution of the paper.
(5) Completion and errors are reported to the trainee as snapshots of missteps.}	
\label{cate}
\end{center}
\end{figure}

\section{Methods}
\subsection{TIPS-author: a surgical simulation creation environment}
\label{sec:author}

%
TIPS-author enables a surgeon-educator (SE) to  improve the specialization, variety and relevance
of laparoscopic VR-training.
TIPS-author  
provides a surgeon-educator (SE)  with an
open source environment to
create and customize hands-on, interactive force-feedback laparoscopy training units. 
Based on the SE's listing, TIPS-author extracts a set of compatible, computationally efficient simlets from a database, and   
generates step-by-step instructions for the web-based TIPS-trainee interface, quizzes and, as specified in Section \secref{sec:safety},  the monitoring of surgical errors. 
The database is generated by a scenario design cycle \cite{Dindar:2016:MMVR} that separates the roles of author, developer and artist:
\begin{itemize}
\item[-]
at the developer level,
numerical simulation routines are selected or adjusted and
admissible parameter ranges are determined;
\item[-]
at the artist level, 
geometric models with collision and physics parameters (simlets) are created; and
\item[-]
at the surgeon-author level, simlets are combined and the workspace and view determined via a interactive webGL interface.
\end{itemize}
Key to this approach is that it can  leverage 
the vibrant professional open-source geometric modeling software \blender\ (http://blender.org,
\cite{blender}) 
and
a library of numerical, geometric and visual algorithms for soft-tissue simulation, \sofa\ (http://www.sofa-framework.org,
\cite{sofaloc,SOFA07}). 
\blender\ and \sofa\ are combined
 to generate hex-meshes on the fly for carvable fatty tissue and thick-shell models of the stomach. 
Anatomy and physics are customized via a menu added to \blender\ that annotates the .xml files sent to \sofa\ for simulation.
Combination and customization 
are  possible 
due to a small disciplined choice of representation primitives.



\subsection{Adding safety monitoring and feedback to TIPS}
\label{sec:safety}
Assessment, evaluation and feedback are critical components in the training of novice surgeons and obeying safety rules is paramount when executing complex sequences of maneuvers.
Physician-training is an experiential process. That is, learners acquire skill by engaging in supervised patient care. All US physicians-in-training, including surgical trainees, must demonstrate competency across a range of knowledge, skills and attitudes prior to graduation \cite{Yardley:2012,Holmboe:2015}. 
Assessing, evaluating and providing critical feedback and instruction
in the workplace is time intensive and stressful. 
And it requires an experienced surgeon's active participation and  expert judgment to provide safe and effective patient care and a quality learning experience.
To ensure that the assessment and evaluation of surgical trainees is reliable and valid 
many training programs employ peer-reviewed evaluative tools such as the 
objective structured assessment of technical skills
(OSATS) for workplace-based assessment \cite{osats:}.

Assessment is also central, but arguably less stressful, in popular computer games where simple counters monitor  progress.  
Psychometric games claim to measure mental agility, attention,  cognitive speed, spatial aptitude and numerical processing ability. 
Increasingly, educational video games incorporate stealth assessment, ubiquitous, unobtrusive, and real-time assessments that intersect play, learning, and assessment. Stealth assessments measure knowledge and skill, then provides learning supports, feedback, instructions, or adapts challenges in the learning environment (e.g., difficulty) to students' proficiency level, maximizing their learning \cite{Shute:2011}. 


Existing VR simulators typically report time to completion, task specific data
such as the number of staples used, and other general counters. 
TIPS's incorporation of SE-established
safety criteria makes cases more relevant for the specific procedure -- but this approach also implies that the \emph{criteria cannot be hard-coded in the simulator ahead of time}. 
Consultation with the clinical experts identified general classes of training errors (i--vi):
\begin{itemize}
\item[i] Incising or cauterizing at the wrong location.
\item[ii] Injuring a nerve by applying too much force (pressure or over-stretching).
\item[iii] Leaving foreign objects in the patient’s body (clips, tools).
\item[iv] Applying surgical clips incorrectly.
\item[v] Removing the wrong (part of) an organ.
\item[vi] Suturing at the wrong location.
\end{itemize}

\begin{figure}[h]
\def \wid{0.45\textwidth}
\def \wwid{0.95\textwidth}
\begin{subfigure}{\wid}
\centering  
\begin{overpic}[trim =  150 100 50 110,clip,scale=.25,tics=10,width=\linewidth]
{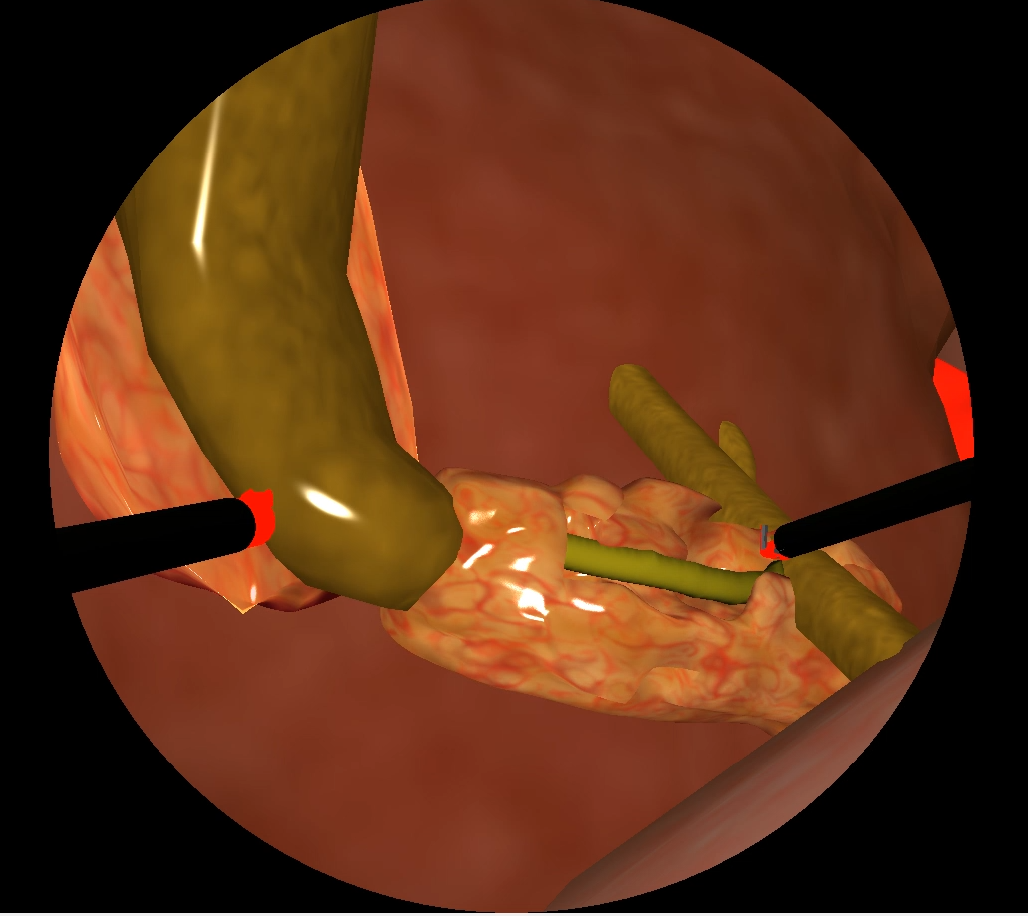}
\put (73,40) {\textcolor{white}{$\Downarrow$}}
\end{overpic}
\caption{Wrong incision on common bile duct.}
\end{subfigure}
\begin{subfigure}{\wid}
\centering    
\begin{overpic}[trim = 150 100 50 100,clip,scale=.25,tics=10,width=\linewidth]
{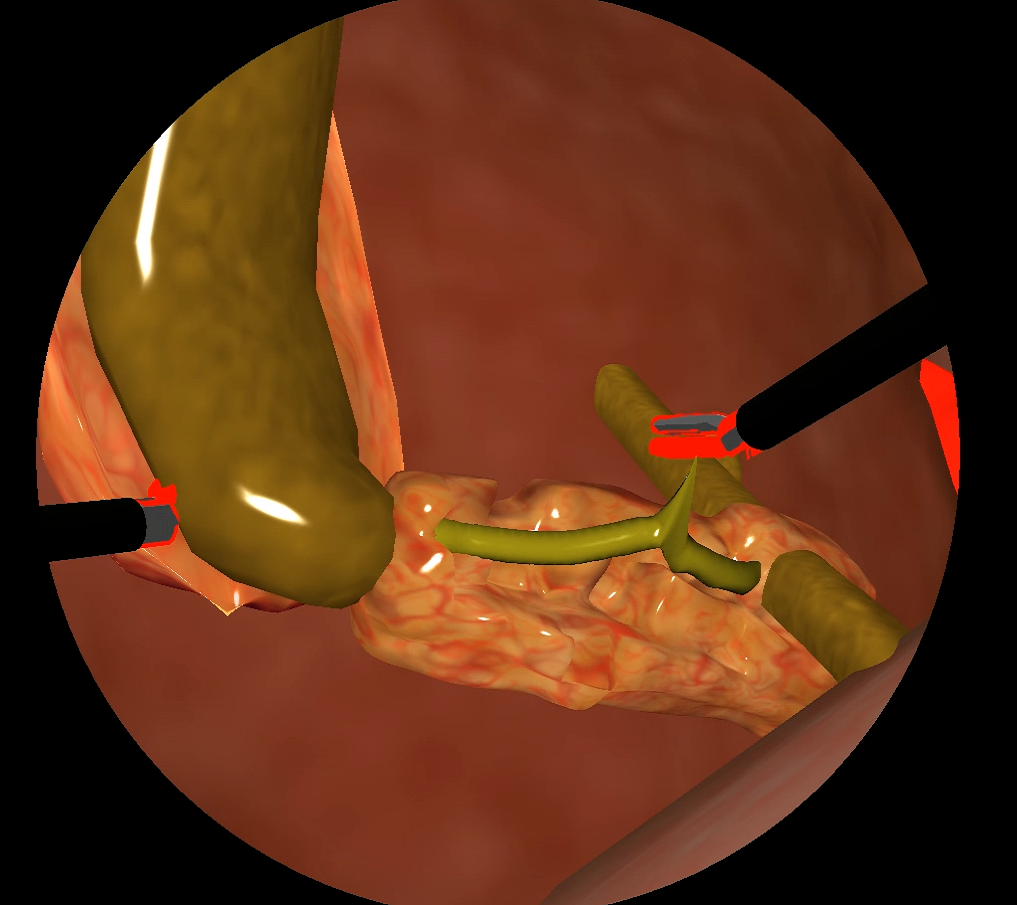}
\put (65,50) {\textcolor{white}{$\Downarrow$}}
\end{overpic}
\caption{Overstretching the cystic duct.}
\end{subfigure}
\begin{subfigure}{\wid}
\centering  
\begin{overpic}[trim =  100 100 100 100,clip, scale=.25,tics=10,width=\linewidth]
{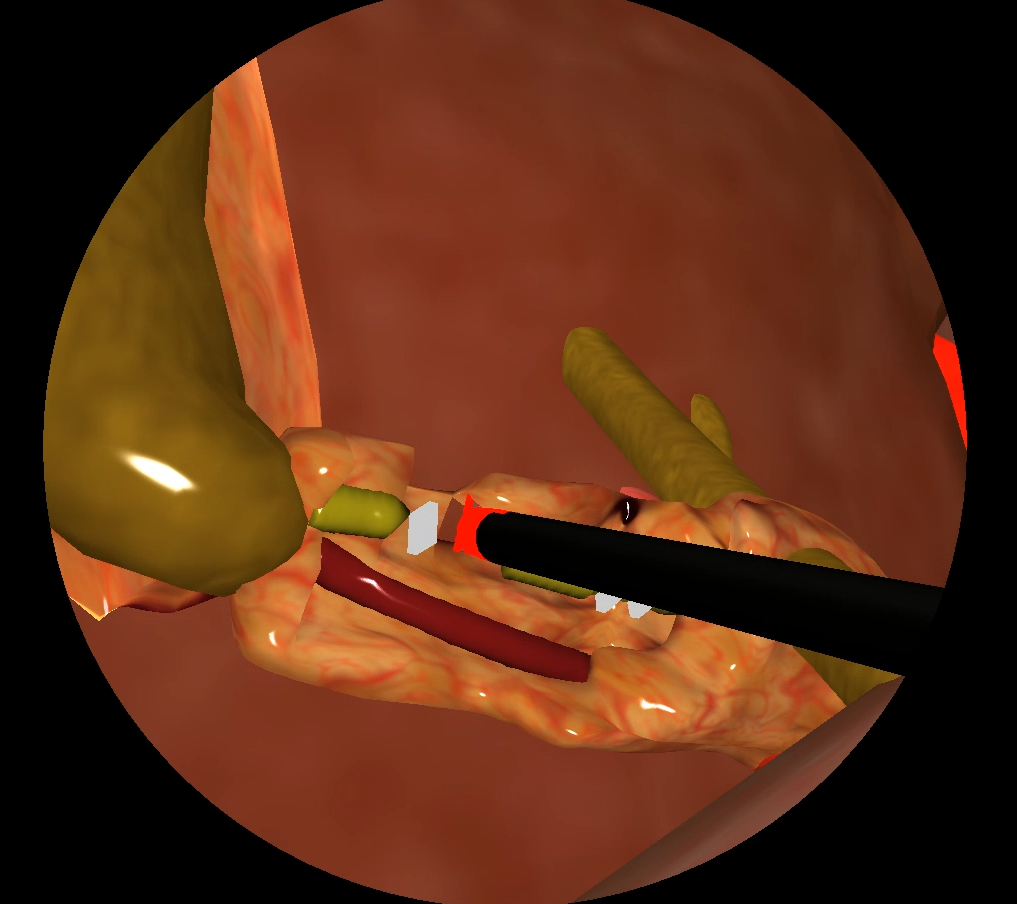}
\put (45,45) {\textcolor{white}{$\Downarrow$}}
\end{overpic}
\caption{Clip drops to abdominal wall due to bile duct cut at the wrong location. }
\end{subfigure}
\hskip0.96cm
\begin{subfigure}{\wid}
\centering  \begin{overpic}[trim =  100 50 100 150,clip,scale=.25,tics=10,width=\linewidth]
{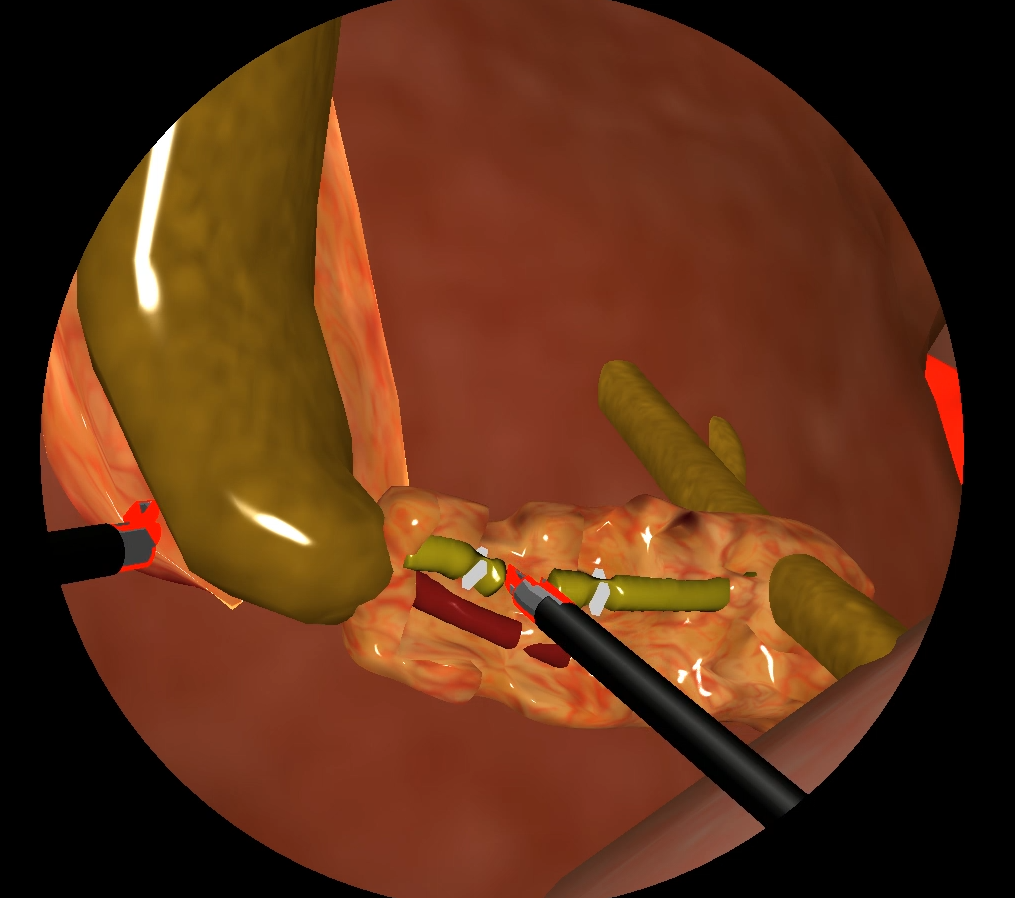}
\put (50,45) {\textcolor{white}{$\Downarrow$}}
\end{overpic}      
\caption{Bile leak due to the lack of vascular clips on the left side.}
\end{subfigure}\hfill   
\caption{Four types of common surgical errors in laparoscopic cholecystectomy reported by TIPS-simulator. For immediate feedback the tool tip becomes red (see $\Downarrow$) 
and the scene briefly freezes
\vid.
}
\label{fig:error}
\end{figure}
These surgical errors can be abstracted as: distance to anatomy, force exerted, location and number of surgical safety clips, and incomplete execution.
Initialized by the `safety' entry in TIPS-author,
our solution is to have TIPS-author parse these safety criteria and  append the corresponding safety tags to these simlets upon export to TIPS-Simulator.
TIPS-simulator monitors these data streams and reports violations both directly and as a sequence of screen-shot images labeled by error types. \figref{fig:error} shows  screen-shots of four common surgical errors (corresponding to type i - iv) during laparoscopic cholecystectomy.

In more detail, error class (i) is monitored by TIPS-simulator as a collision event with an offset distance between the tool listed in the TIPS-authors tuple of the task and an organ listed in the safety entry.

For example, for cholecystectomy, for the step `Explore the triangle of Calot' (see \figref{fig:error}a) the task tuple reads:
\begin{itemize}
\item[]
dissect Fatty tissue over the cystic ductus and cystic artery\\
using Curved Maryland Dissector not too close to Common bile duct.
\end{itemize}

\begin{wrapfigure}{r}{0.35\textwidth}
\centering  \includegraphics[width=\linewidth]{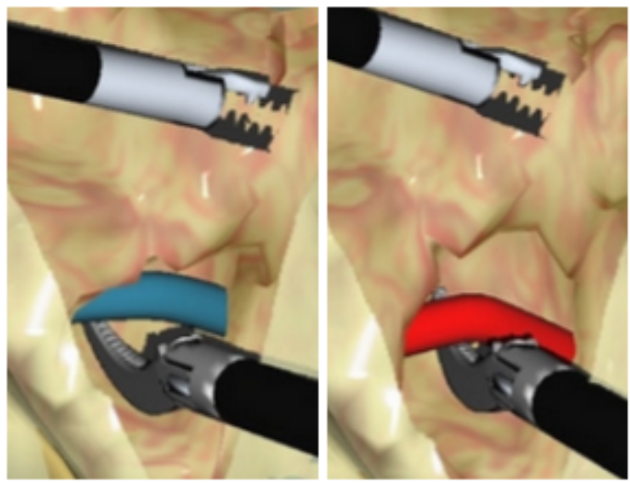}
\caption{Excessive force indicated by change of vessel to red (flashing) color.}
\label{fig:force}
\end{wrapfigure}
Here `dissect' is the action, `Fatty tissue over the cystic ductus and cystic artery' is the anatomy (specifying a simlet), `Curved Maryland Dissector' specifies the laparoscopic tool,  `not too close' indicates an error of type (i) and `Common bile duct' is an entry in the anatomy database that requires monitoring.
TIPS-simulator then monitors distance between the cauterizing tool and the common bile duct. Distance below the offset triggers and registers the error.

Type (ii) errors are monitored in terms of force feedback returned to the haptic devices. This safety threshold is customized for veins or arteries with different physical properties.
Type (iii) and (iv) errors are detected by a map tracking the vector of deployed clips on each clip-able object to  monitor not only the number but also placement of clips. For example, to prevent bleeding or leaking, two clips should be applied on the part of the duct or vein that remains inside the body. 
Type (v) errors are indirectly caught since they terminate the simulation
without generating an `achievement' entry in the final visual report and type (vi) errors are caught by initializing suturable regions on the object, say the fundus of the gaster during Nissen fundoplication.

Errors (i - iv) alert the trainee  by a red flashing (\figref{fig:error}a) instrument tip.
A corresponding screen-shot is saved for later named by time, error type, and error values. 

\subsection{Summary feedback as a series of snapshots}
Once the procedure completes, typically when the cancerous organ is retrieved via the surgical pouch, all screen-shots of errors (and small ones for task completions) are displayed to the trainee as a feedback report. This serves as starting point for a discussion with the instructor. 
Proficiency with respect to  the training module is equivalent to repeated performance without errors and a complete list of achievements. 
The final achievement is generically checked by asserting that the cancerous body is free from the remaining organs and tissues. Similarly, clip placement requires freeing the vessel and testing that two clips remain within the body while a third clip ensures integrity of the tissue to be removed. Such authored criteria provide more valuable feedback than time taken or number of clips deployed.

\begin{wrapfigure}{R}{0.4\textwidth}
\begin{center}
\vskip-1.0cm
\includegraphics[width=\linewidth]{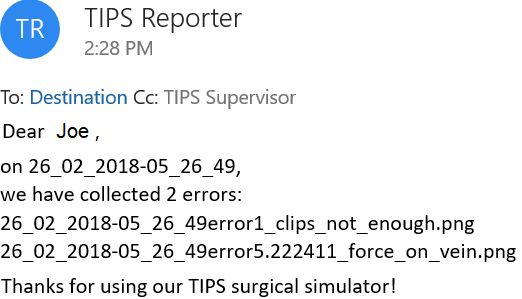}
\caption{Email reports to trainee (and the instructor) upon completion.}
\vskip-1.0cm
\label{fig:email}
\end{center}
\end{wrapfigure}
Additionally the unique directory of screen-shots and the filenames are reported to the trainee by e-mail  and, optionally, to an account set up for the  instructor (see \figref{fig:email})
to document training progress and decide whether the pattern and number of errors
requires intervention and what errors should be discussed.

In summary, faced with the complexity of supporting procedure-specific proficiency assessment, we categorized laparoscopic safety violations into several generic classes. This enables a simple but effective, implemented and tested strategy to use TIPS-author safety entries to initialize, monitor and report error events; and to create a record of progress towards proficiency. 


\section{Results}
\subsection{Evaluation of TIPS training and feedback}

TIPS was demonstrated and experienced by a broad range of medical professionals at the
the American College of Surgeons Clinical 
Congress 2019 (ACS) and the Academic Surgical Congress 2020 (ASC).
Besides testing the technology `in the wild', the venues allowed the team  to conduct a survey of TIPS.

\begin{figure}[h]
\begin{subfigure}{0.45\textwidth}
\centering        
\includegraphics[width=0.96\textwidth]{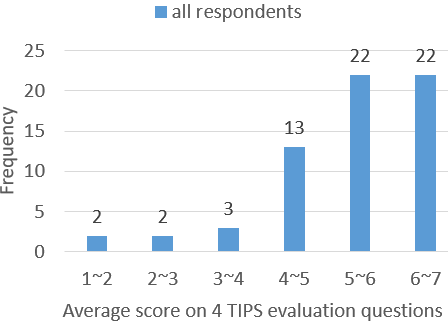}
\caption{cumulative}
\end{subfigure}
\begin{subfigure}{0.55\textwidth}
\centering        
\includegraphics[width=0.96\textwidth]{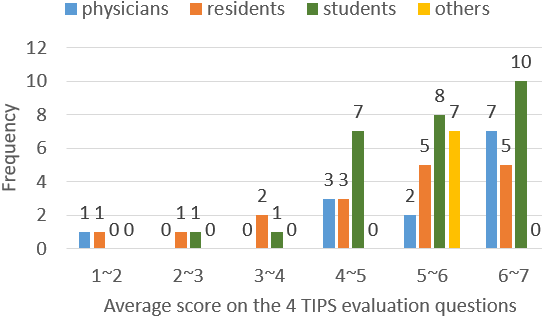}
\caption{by group}
\end{subfigure}
\caption{Breakdown of the average score on the 4 central TIPS evaluation questions.}
\label{fig:frequency}
\end{figure}

Prior to field testing, face, construct and content validity of SE-authored cholecystectomy and appendectomy TIPS modules had been established by laparoscopic surgeons and residents at the Universities of Florida and Buffalo.
At the congress field tests, after training with the modules, 64 respondents (13 board certified surgeons, 17 medical residents, 27 medical students and 7 other medical professionals) rated TIPS across several usability items
on a Likert scale from 1 to 7  (7 = strongly agree, see \cite{surveyACS}).
The scale resolution was selected as a trade-off between scale complexity and expressiveness.
Table \ref{tab:qs} lists the outcome of the four central questions of usability and \figref{fig:frequency} breaks down the score on these four central questions. (Four other questions established medical seniority, familiarity with virtual trainers, and prior experience with laparoscopy). 
All questions were selected in consultation with SEs at the authors' institutions.

\begin{table}
\begin{tabular}
{p{0.75\textwidth}
p{0.1\textwidth}
p{0.1\textwidth}
}
& mean & standard 
\\
TIPS  ...   & rating  & deviation
\\
\hline
helps understand the force necessary to explore the anatomy   & 5.34 & 1.46
\\
\hline
interface does not distract from the surgical task & 5.02 & 1.52\\
\hline
enhances lap-competency attainment over current methods
& 5.19 & 1.5\\
\hline
is compatible with the current lap training curricula & 5.39 &
1.43 \\
\hline
overall score  & 5.24 & 1.33\\
\end{tabular}
\caption{TIPS with safety rules rated on the four key questions.}
\label{tab:qs}
\end{table}

\subsection{The effect of summative visual feedback via snapshots}
When analyzing the data sets from the two conferences, we noticed agreement
of averages between the 13 ACS respondents and the 51 ASC respondents. The agreement 
was within .2 on all
rating categories except one. 
The only outlier was the statement 
`The TIPS interface helps learners understand the force necessary to explore the anatomy'. Here the rating  improved from 5.04 at ACS to 5.35 at ASC.
The only change applied to the TIPS software after the ACS survey
and before the ASC survey,
was the addition of the visual summary of the achievements and procedural errors as a series of snapshots. 
The immediate feedback by changing the color of a vessel was present in both tests.


\section{Discussion}
TIPS is a novel authoring environment that allows surgeon-educators to build customizable VR lap scenarios.	
The bulk of the survey questions was aimed to  evaluate TIPS as a whole.
Ratings collected from medical professionals at two conference exhibitions indicate the viability of such SE-authored surgical training.
In particular, the high score  for `enhances lap-competency attainment over current methods' speaks to the added value of customized TIPS simulations over available current methods.

We did not set out to measure the impact of automatic visual summative feedback on errors presented as screen snapshots.
In fact, the immediate feedback on SE-authored error measurements, a change of color, was present at both field tests.
It is therefore noteworthy that presenting trainee errors in visual form at the end of the training increased acceptance noticeably.
Indeed, additional informal feedback
from surgeons and trainees to the question 'what feature of TIPS do you recall'
endorsed visual feedback via screen shots as both  meaningful and memorable.








\bibliographystyle{elsarticle-num}
\bibliography{p}

\begin{thebibliography}{10}
\expandafter\ifx\csname url\endcsname\relax
  \def\url#1{\texttt{#1}}\fi
\expandafter\ifx\csname urlprefix\endcsname\relax\def\urlprefix{URL }\fi
\expandafter\ifx\csname href\endcsname\relax
  \def\href#1#2{#2} \def\path#1{#1}\fi

\bibitem{Macario:2010}
A.~Macario, What does one minute of operating room time cost?, J Clin Anesth.
  22~(4) (2010) 233--6.
\newblock \href {https://doi.org/10.1016/j.jclinane.2010.02.003. PMID:
  20522350} {\path{doi:10.1016/j.jclinane.2010.02.003. PMID: 20522350}}.

\bibitem{Zendejas:2013}
Z.~B, Brydges, H.~S.~C. DA., State of the evidence on simulation-based training
  for laparoscopic surgery: a systematic review, Ann Surg. 257~(4) (2013)
  586--593, pMID: 23407298.
\newblock \href {https://doi.org/10.1097/SLA.0b013e318288c40b}
  {\path{doi:10.1097/SLA.0b013e318288c40b}}.

\bibitem{qin2010novel}
J.~Qin, W.-M. Pang, Y.-P. Chui, T.-T. Wong, P.-A. Heng, A novel modeling
  framework for multilayered soft tissue deformation in virtual orthopedic
  surgery, Journal of medical systems 34~(3) (2010) 261--271.

\bibitem{talbot2015surgery}
H.~Talbot, N.~Haouchine, I.~Peterlik, J.~Dequidt, C.~Duriez, H.~Delingette,
  S.~Cotin, Surgery training, planning and guidance using the sofa framework,
  in: Eurographics, 2015.

\bibitem{echegaray2014brain}
G.~Echegaray, I.~Herrera, I.~Aguinaga, C.~Buchart, D.~Borro, A brain surgery
  simulator, IEEE computer Graphics and Applications 34~(3) (2014) 12--18.

\bibitem{arikatla2018high}
V.~S. Arikatla, M.~Tyagi, A.~Enquobahrie, T.~Nguyen, G.~H. Blakey, R.~White,
  B.~Paniagua, High fidelity virtual reality orthognathic surgery simulator,
  in: Medical Imaging 2018: Image-Guided Procedures, Robotic Interventions, and
  Modeling, Vol. 10576, International Society for Optics and Photonics, 2018,
  p. 1057612.

\bibitem{mitchell2015gridiron}
N.~Mitchell, C.~Cutting, E.~Sifakis, Gridiron: an interactive authoring and
  cognitive training foundation for reconstructive plastic surgery procedures,
  ACM Transactions on Graphics (TOG) 34~(4) (2015) 1--12.

\bibitem{cecil2018advanced}
J.~Cecil, A.~Gupta, M.~Pirela-Cruz, An advanced simulator for orthopedic
  surgical training, International journal of computer assisted radiology and
  surgery 13~(2) (2018) 305--319.

\bibitem{Satava2005}
A.~Gallagher, E.~Ritter, H.~Champion, G.~Higgins, M.~Fried, G.~Moses, C.~Smith,
  R.~Satava, Virtual reality simulation for the operating room:
  proficiency-based training as a paradigm shift in surgical skills training,
  An. Surgery 241~(2) (2005) 364--372.

\bibitem{Gurusamy2009}
L.~P. K.S.~Gurusamy, R.~Aggarwal, B.~Davidson, Virtual reality training for
  surgical trainees in laparoscopic surgery, The Cochrane Database of
  Systematic Reviews 1 (2009).

\bibitem{Yeo:2010:MMVR}
Y.-I. Yeo, S.~Dindar, G.~Sarosi, J.~Peters, Enabling surgeons to create
  simulation-based teaching modules, in: {M}edicine {M}eets {V}irtual {R}eality
  (MMVR) Feb 8-12 2011, Long Beach,CA, Studies in Health Technology and
  Informatics (SHTI), IOS Press, Amsterdam, 2011, pp. 1--6.

\bibitem{Dindar:2016:MMVR}
S.~Dindar, T.~Nguyen, J.~Peters, Towards surgeon-authored {VR} training: the
  scene-development cycle, in: Proceedings of {M}edicine {M}eets {V}irtual
  {R}eality (MMVR) April 9-12 2016, L.A. ,CA, Studies in Health Technology and
  Informatics (SHTI), IOS Press, Amsterdam, 2016, pp. 1--6.

\bibitem{Sarov:2018:AIS}
M.~Sarov, R.~Gao, J.~Youngquist, G.~Sarosi, S.~Kurenov, J.~Peters, An authoring
  interface for surgeon-authored {VR} training, International Journal of
  Computer Assisted Radiology and Surgery 13 (2018) 1--(1--4)--273.

\bibitem{Lesch:2019:VR:}
H.~Lesch, E.~Johnson, J.~Cendan, J.~Peters, {VR} simulation leads to enhanced
  procedural confidence for surgical trainees, Journal of Surgical Education
  77~(1) (2019) 213--218, nIHMS1680076.

\bibitem{wise}
N.~G.~S. of~Medicine, A.~C. of~Surgeons~(ACS), A.~of~Surgical Education~(ASE),
  \href{http:/www.aquifer.org}{Wise-md}.
\newline\urlprefix\url{http:/www.aquifer.org}

\bibitem{McGrane:2016}
M.~B, Belton, Powell, I.~J., The relationship between fundamental movement
  skill proficiency and physical self-confidence among adolescents, J Sports
  Sci 35~(17) (2017) 1709--1714, pMID: 28282760.
\newblock \href {https://doi.org/10.1080/02640414.2016.1235280}
  {\path{doi:10.1080/02640414.2016.1235280}}.

\bibitem{blender}
B.~Foundation, Blender, \url{http://blender.org/} (2015).

\bibitem{sofaloc}
INRIA, Simulation open framework architecture,
  \url{http://www.sofa-framework.org} (2015).

\bibitem{SOFA07}
J.~Allard, S.~Cotin, F.~Faure, P.-J. Bensoussan, F.~Poyer, C.~Duriez,
  H.~Delingette, L.~Grisoni, Sofa \- an open source framework for medical
  simulation, in: Medicine Meets Virtual Reality (MMVR'15), Long Beach, USA,
  2007.

\bibitem{Yardley:2012}
S.~Yardley, P.~W. Teunissen, T.~Dornan, Experiential learning: transforming
  theory into practice, Medical teacher 34~(2) (2012) 161--164.

\bibitem{Holmboe:2015}
E.~S. Holmboe, Realizing the promise of competency-based medical education,
  Academic Medicine 90~(4) (2015) 411--413.

\bibitem{osats:}
J.~Martin, G.~Regehr, R.~Reznick, H.~MacRae, J.~Murnaghan, C.~Hutchison,
  M.~Brown, Objective structured assessment of technical skill (osats) for
  surgical residents, Br J Surg. (1997) 273--8.

\bibitem{Shute:2011}
V.~J. Shute, Stealth assessment in computer-based games to support learning,
  Computer games and instruction 55~(2) (2011) 503--524.

\bibitem{surveyACS}
\href{https://ufl.qualtrics.com/jfe/form/SV\_b0WFv9hk3vl44fj}{Survey for 2019
  {A}merican {C}ollege of {S}urgeons {C}linical {C}ongress and the 2020
  {A}cademic {S}urgical {C}ongress}.
\newline\urlprefix\url{https://ufl.qualtrics.com/jfe/form/SV\_b0WFv9hk3vl44fj}

\end{thebibliography}
\end{document}